\documentclass[12pt, reqno, twoside, letterpaper]{amsart}
\usepackage{amsmath,amssymb,amsbsy,amsfonts,amsthm,latexsym,
            amsopn,amstext,amsxtra,euscript,amscd,stmaryrd,mathrsfs,
            cite,array}

\newtheorem{X}{X}[section]
\newtheorem{theorem}[X]{Theorem}  

\newtheorem{corollary}[X]{Corollary}

\newtheorem{proposition}[X]{Proposition}

\makeatletter
\renewcommand{\pmod}[1]{\allowbreak\mkern7mu({\operator@font mod}\,\,#1)}
\makeatother

\newcommand{\order}{\asymp}      
\renewcommand{\(}{\left(}
\renewcommand{\)}{\right)}

\newcommand{\be}{\begin{equation}}
\newcommand{\ee}{\end{equation}}
\newcommand{\dalign}[1]{\[\begin{aligned} #1 \end{aligned}\]}

\newcommand{\eps}{\varepsilon}

\newcommand{\cN}{\EuScript{N}}

\title{Large prime gaps and progressions with few primes}

\author[K.~Ford]{Kevin Ford}

\address{Department of Mathematics,
         University of Illinois, 1409 West Green St,
         Urbana, IL 61801, USA.}

\email{ford126@illinois.edu}

\date{\today}

\begin{document}

\begin{abstract}
We show that the existence of 
arithmetic progressions with few primes,
with a quantitative bound on ``few'',
implies  the existence of larger
gaps between primes less than $x$ than is currently known unconditionally.  In particular, 
we derive this conclusion if there are certain types of
exceptional zeros of Dirichlet $L$-functions.
\end{abstract}

\thanks{MSC Primary: 11N05, 11N13; Secondary 11M20.}


\thanks{\textbf{Acknowledgments:} The author was supported by National Science Foundation grant DMS-1802139.}
\maketitle



{\section{Introduction}
\label{sec:intro}}

Estimation of the largest gap, $G(x)$, between consecutive primes less than $x$ is a classical problem, and 
the best bounds on $G(x)$ are comparatively weak.
The strongest unconditional lower bound on $G(x)$
is due to Ford, Green, Konyagin, Maynard and Tao~\cite{FGKMT}, who
have shown that
\be\label{FGKMT}
G(x)\gg\frac{\log x\,\log_2x\,\log_4x}{\log_3x},
\ee
for sufficiently large $x$, with $\log_k x$ the $k$-fold iterated natural logarithm of $x$, whereas the best unconditional upper bound is
\be\label{BHP}
G(x)\ll x^{0.525},
\ee
a result due to Baker, Harman and Pintz~\cite{BHP}.  Assuming the Riemann Hypothesis, Cram\'er~\cite{Cramer1} showed that
\[
G(x) \ll x^{1/2} \log x.
\]
The huge distance between the lower bound \eqref{FGKMT} and upper bound
\eqref{BHP} testifies to our ignorance about gaps between primes.
Cram\'er \cite{Cra36} introduced a probabilistic model for primes
and used it to  conjecture that $\limsup_{x\to \infty} G(x)/\log^2 x \ge 1$;
later, Shanks \cite{Shanks} conjectured that $G(x)\sim \log^2 x$ based on
a similar model.
Granville \cite{Gra95} modified Cram\'er's model
and, based on analysis of the large gaps in the model, conjectured that $G(x)\ge (1+o(1))2e^{-\gamma} (\log x)^2$.  The author, together with William Banks and Terence Tao \cite{BFT}, 
has created another model of primes $\le x$,
the largest gap in the model set depending
on an extremal property of the interval sieve.
In particular, the existence of a certain sequence of ``exceptional zeros''
of Dirichlet $L$-functions (defined below) implies that 
the largest gap in the model set grows faster than
any constant multiple of $(\log x)^2$, and suggests that 
the same bound holds for $G(x)$.
 In this paper, we show that the existence of exceptional zeros
of a certain type implies a lower bound for $G(x)$
which is larger than the right side of \eqref{FGKMT}.
We do not utilize probabilistic models of primes, but instead we argue directly.
More generally, we derive a similar conclusion whenever
there are arithmetic progressions containing few primes.
We denote $\pi(x;q,b)$ the number of primes $p\le x$
satisfying $p\equiv b\pmod{q}$.
The prime number theorem for arithmetic progressions implies that
\[
\pi(x;q,b) \sim \frac{\pi(x)}{\phi(q)}
\]
for any \emph{fixed} $q$, where $\phi$ is Euler's totient function
and $\pi(x)$ denotes the number of primes $p\le x$.  It is a central problem to prove bounds on
$\pi(x;q,b)$ which are uniform in $q$, but 
the best known results are only uniform for 
$q\le (\log x)^{O(1)}$; see
\cite{Dav} for the classical theory.

All of the methods used to prove lower bounds on $G(x)$
utilize a simple connection between $G(x)$ and
Jacobsthal's function $J(u)$, the maximum gap between integers
having no prime factor $p\le u$.  A simple argument based on 
the prime number theorem and the
 Chinese Remainder Theorem implies that
\be\label{Jacobsthal}
G(x) \ge J((1/2)\log x)
\ee
if $x$ is sufficiently large.
The best bounds known today for $J(u)$ are
\[
u (\log u) \frac{\log_3 u}{\log_2 u} \ll J(u) \ll u^2,
\]
the lower bound proved in \cite{FGKMT} and the upper bound
due to Iwaniec \cite{I71}.

\medskip

\begin{theorem}\label{thm: AP_gaps}
Suppose that $x$ is large, $x>q>b>0$ and
 $\pi(x;q,b) \le  \frac{\delta x}{\phi(q)}$
with $0\le \delta\le 1$.  Then
\[
G(e^{2u}) \ge J(u) \ge \frac{x-b}{q}.
\]
where $u$ is the smallest integer satisfying $u>2\sqrt{x}$
and $\frac{u}{\log u} \ge \frac{10\delta x}{q}$. 
\end{theorem}

An immediate corollary gives a lower bound on $G(x)$
assuming a lower bound on $L(q,b)$, the least prime
in the progression $b\mod q$.  We take $\delta=0$ and $u=\lceil 2\sqrt{x} \rceil$.

\begin{corollary}\label{cor:Linnik}
Suppose that $L(q,b) > x$.  Then $G(e^{4\sqrt{x}})\ge \frac{x-b}{q}$.
\end{corollary}

Theorem \ref{thm: AP_gaps} is a partial converse
to a theorem of Pomerance \cite[Theorem 1]{Pom},
which provides a lower bound on $\max_{(b,q)=1} L(q,b)$
given a lower bound on the maximal gap between numbers
coprime to $m$, where $(m,q)=1$ and  $m\le q^{1-o(1)}$.

Linnik's theorem \cite{Linnik} states that $L(q,b) \ll q^L$ for some
constant $L$; the best quantitative result of this kind
is
due to Xylouris \cite{Xylouris}, who showed that 
the bound holds with $L=5.18$.
Assuming the Extended Riemann Hypothesis (ERH) for
Dirichlet $L$-functions, we obtain  a stronger bound 
$L(q,b) \ll_\eps q^{2+\eps}$ for every $\eps>0$.
If, for some $c>2$ there are infinitely many pairs $(q,b)$
with $L(q,b)\ge q^c$ (a violation of ERH), then Corollary \ref{cor:Linnik}
implies that
\[
\limsup_{X\to \infty} \frac{G(X)}{(\log X)^{2-\frac{2}{c}}} > 0.
\]
It is, however, conjectured that $L(q,b) \ll \phi(q) \log^2 q$;
see \cite{LPS} for a precise version of this conjecture
and for the best known lower bounds on $\max_{(b,q)=1} L(q,b)$.

We may also exceed the bound in \eqref{FGKMT} under
the assumption that exceptional zeros of Dirichlet
$L$-functions exist.
Roughly speaking, an exceptional zero 
of $L(s,\chi)$ is a zero which is real
and very close to 1.  As such, their existence violates
ERH for $L(s,\chi)$.
Classical results (see \cite[\S 14]{Dav}) imply that
if $c_0>0$ is small enough, and $q\ge 3$, then there is at most one character $\chi$ modulo $q$ for which $L(s,\chi)$
has a zero in the region $$\{ \sigma+it\in \mathbb{C} : \sigma \ge 1-c_0/\log(qt) \},$$ and moreover
the character is real and the zero is real.
We shall refer to such zeros as ``exceptional zeros''
with respect to $c_0$.
Moreover, by reducing $c_0$ if necessary, 
it is known that moduli $q$ for which an exceptional zero exists
are very rare.

  Siegel's theorem
\cite[Sec. 21]{Dav} implies that 
\be\label{Siegel}
\log \frac{1}{1-\beta_q} = o(\log q) \quad (q\to \infty),
\ee
for (hypothetical) exceptional zeros $\beta_q$,
 although we cannot say any rate at which this
occurs (the bound is \emph{ineffective}). 
The exceptional zeros are also know as
Siegel zeros or Landau-Siegel zeros in the literature.
Their existence implies a great irregularity
in the distribution of primes modulo $q$, given by
Gallagher's Prime Number Theorem \cite{Gal70}.
Here we record an immediate corollary.

\begin{proposition}[Gallagher]
\label{Gallagher}
For some absolute constant $B>1$, we have the following.
Suppose that $\chi$ is a real character with conductor $q$ and
$L(1-\delta,\chi)=0$ for some $0<\delta<1$.
Then, for all $b$ with  $\chi_q(b)=1$ and all $x\ge q^B$, we have
\[
\pi(x;q,b) \ll \frac{\delta x}{\phi(q)}.
\]
\end{proposition}

One can leverage this irregularity to prove
\emph{regularity} results about primes that are
out of reach otherwise, 
the most spectacular application being Heath-Brown's
\cite{HB}
deduction of the twin prime conjecture from
the existence of exceptional zeros (for an appropriate $c_0$).
See Iwaniec's survey article \cite{Iw02 } 
 for background on attempts to prove the non-existence of
 exceptional zeros and discussion
about other applications of their existence.
There are also a variety of problems where one
argues in different ways depending on whether or not
exceptional zeros  exist, a principal example being
Linnik's Theorem on primes in arithmetic progressions
(see, e.g., \cite[Ch. 24]{FI}).

Apply Proposition \ref{Gallagher} with $x=q^B$. 
Recalling \ref{Siegel}, we see that the quantity $u$
in Theorem \ref{thm: AP_gaps} satisfies
\[
u \order \frac{\delta x \log x}{q}
\]
and consequently that $\log u \order \log q$.
We conclude that

\begin{theorem}\label{thm:main}
Suppose that $\chi$ is a real character with conductor $q$ and that
$L(1-\delta,\chi)=0$ for some $0<\delta<1$.
Then
\be\label{G(x)_Siegel}
G(e^{2u}) \gg  \frac{u}{\delta \log u},
\ee
for some $u$ satisfying $\log u \order \log q$.
\end{theorem}

For example, if $k$ is fixed and there exist infinitely many
exceptiona zeros $\delta=\delta_q$ satisfying
 $\delta_q \le (\log q)^{-k}$, we see that
there is an unbounded set of $X$ for which
\[
G(X) \gg_k (\log X)(\log_2 X)^{k-1}.
\]
this improves upon \eqref{FGKMT} for $k\ge 2$.
Similarly, if there is an infinite set of $q$ satisfying
$\delta=\delta_q=q^{-\eps(q)}$, where $\eps(q)\to 0$ very slowly,
then for an unbounded set of $X$,
\[
G(X) > X^{1+\delta(X)}
\]
with $\delta(X)\to 0 $ very slowly.


{\section{Proof of Theorem \ref{thm: AP_gaps}}}

Let $u$ be as in the theorem, and  let
\be\label{yz}
y = \frac{x-b}{q}.
\ee
Too show that $J(u) \ge y$, it suffices to
find residue classes $a_p\mod p$, one for each prime $p\le u$, which 
together cover $[0,y]$. 
For each prime $p\le u/2$ with $p\nmid q$, define
 $a_p$ by
\[
q a_p + b \equiv 0\pmod{p}.
\]
Recall that $(b,q)=1$.
In this way, if $0\le n\le y$ and $n\not\equiv a_p\pmod{p}$ 
for all such $p$, then $m=qn+b$ has no prime factor $\le u/2$.
Also, $x=qy+b<(u/2)^2$ by hypothesis,
and thus $m$ is prime. Let
\[
\cN = \{ 0\le n\le y : n \not\equiv a_p \pmod{p}, \, \forall p\le u/2 \text{ with }p\nmid q\}.
\]
It follows from the hypothesis of the theorem that
\dalign{
|\cN| &\le \pi(qy+b;q,b)=\pi(x;q,b) \le 
\frac{\delta x}{\phi(q)}.
}
Next, we choose residue classes $a_p$ for primes $p|q$ with $p\le u/2$ using a greedy algorithm, successively selecting 
for each $p$ a residue class $a_p\mod p$ which covers at least a proportion $1/p$ of the elements remaining uncovered.
As $u> 2\sqrt{x} > 2\sqrt{q}$, there is
at most one prime $p|q$ satisfying $p>u/2$.
Letting $\cN'$ denote the set of
$n\in [0,y]$ not covered by $\{a_p\mod p: p\le u/2\}$, we have
\[
|\cN'| \le |\cN| \, \prod_{p|q,p\le u/2}\(1-\frac{1}{p} \)
\le 2 |\cN| \frac{\phi(q)}{q} \le \frac{2\delta x}{q}.
\]
By hypothesis,
\[
|\cN'| \le \frac{u}{5\log u},
\]
which, by the prime number theorem,
is less than the number of primes
in $(u/2,u]$ for $u$ large enough (as $u>\sqrt{x}$, this happens if $x$ is large enough). Thus,
 we may associate each number  $n\in \cN'$ with a distinct prime
in $p_n\in(u/2,u)$.  Choosing $a_{p_n}\equiv n\pmod{p_n}$
   for each $n\in \cN'$ then ensures that
$\{a_p\mod p: p\le u\}$ covers all of $[0,y]$, as desired.
\qed

\textbf{Remark.}  We have made no
use in the proof of estimates for numbers lacking large prime factors, a common feature in unconditional lower bounds on
$G(x)$.  There does not seem to be any advantage to this in
our argument.




\end{document}